\documentclass[final,12pt]{elsarticle}

\usepackage{amssymb}
\usepackage{amsmath}
\usepackage{algorithm}  
\usepackage{booktabs}
\usepackage{enumitem}   
\usepackage{amsthm}
\usepackage{algorithm}
\usepackage{algpseudocode}

\journal{Journal of Computational Physics}

\begin{document}

\begin{frontmatter}



\title{B-spline periodization of Fourier pseudo-spectral method for non-periodic problems
}


\author[THU]{Dongan Li} 
\author[MPCDF]{Mou Lin}
\author[THU]{Shunxiang Cao}
\author[THU]{Shengli Chen\corref{cor1}}
\affiliation[THU]{organization={Institute for Ocean Engineering, Shenzhen International Graduate School, Tsinghua University},
            addressline={Lishui Road 2279}, 
            city={Shenzhen},
            postcode={518055},
            country={China}}

\affiliation[MPCDF]{organization={Max Planck Computing and Data Facility},
            addressline={Gießenbachstraße 2}, 
            city={Garching},
            postcode={85748}, 
            country={Germany}}
\cortext[cor1]{Corresponding author: shenglichen@sz.tsinghua.edu.cn}
\begin{abstract}
Spectral methods are renowned for their high accuracy and efficiency in solving partial differential equations. The Fourier pseudo-spectral method is limited to periodic domains and suffers from Gibbs oscillations in non-periodic problems.  The Chebyshev method mitigates this issue but requires edge-clustered grids, which does not match the characteristics of many physical problems. To overcome these restrictions, we propose a B-spline-periodized Fourier (BSPF) method that extends to non-periodic problems while retaining spectral-like accuracy and efficiency. The method combines a B-spline approximation with a Fourier-based residual correction. The B-spline component enforces the smooth matching of boundary values and derivatives, while the periodic residual is efficiently treated by Fourier differentiation/integration. This construction preserves spectral convergence within the domain and algebraic convergence at the boundaries. Numerical tests on differentiation and integration confirm the accuracy of the BSPF method superior to Chebyshev and finite-difference schemes for interior-oscillatory data. Analytical mapping further extends BSPF to non-uniform meshes, which enables selective grid refinement in regions of sharp variation. Applications of the BSPF method to the one-dimensional Burgers' equation and two-dimensional shallow water equations demonstrate accurate resolution of sharp gradients and nonlinear wave propagation, proving it as a flexible and efficient framework for solving non-periodic PDEs with high-order accuracy.
\end{abstract}

\begin{graphicalabstract}
\includegraphics[width=1\linewidth]{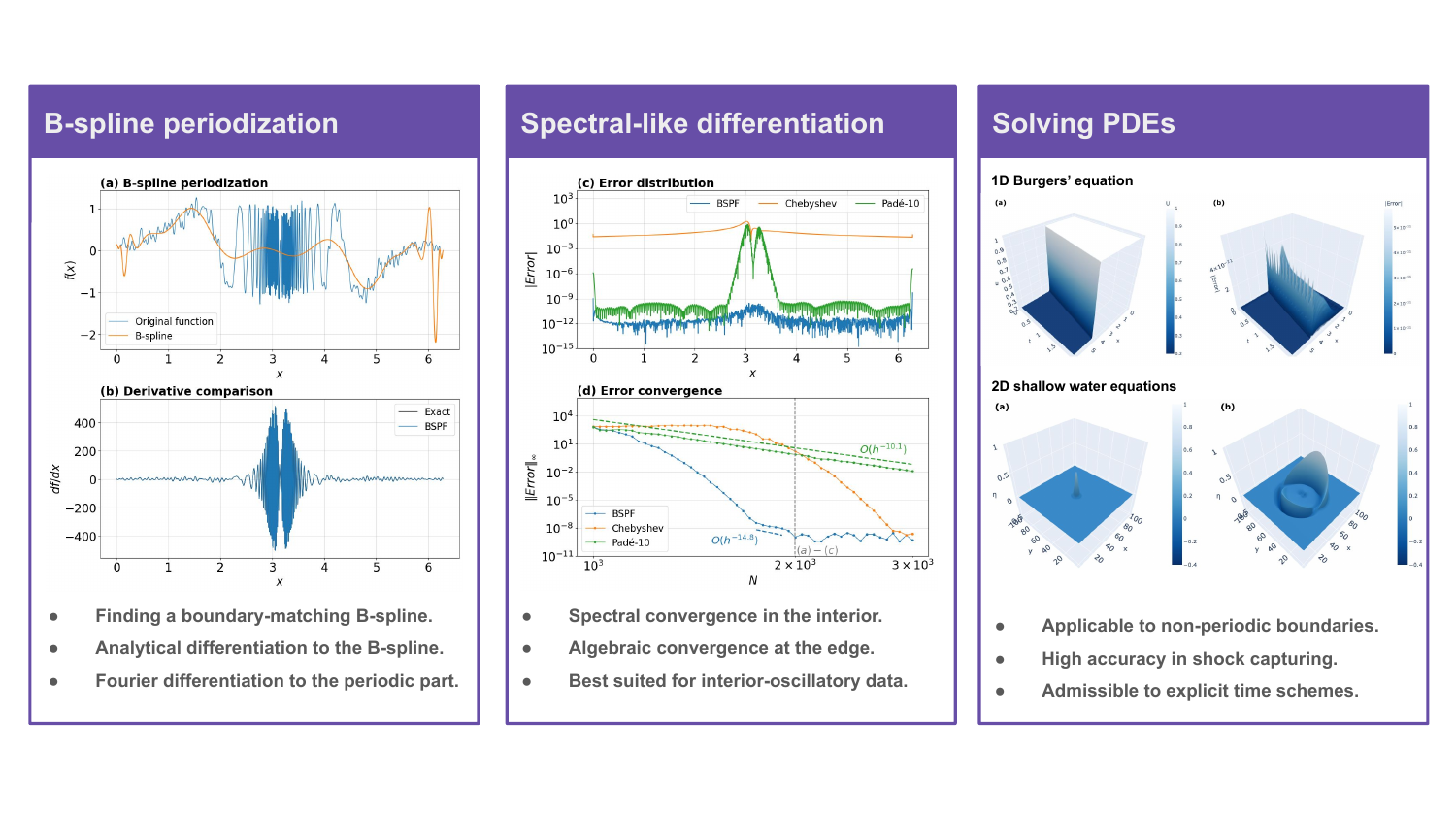}
\end{graphicalabstract}

\begin{highlights}
\item A B-spline-periodized Fourier (BSPF) framework is introduced to achieve boundary periodization with high-order smoothness for non-periodic problems.
\item Boundary conditions are matched by regularized B-spline constraints, and periodic residuals are treated by FFT, so near-spectral accuracy and the high $O(N\log(N))$ efficiency are retained in the BSPF method.
\item The BSPF method achieves superior accuracy over the Chebyshev and finite difference methods for differentiation and integration on interior-oscillatory data, with similar computational cost to the Chebyshev method.
\item Sharp shock waves in the solutions of 1D Burgers' and 2D shallow-water equations are well-resolved, and efficient explicit time stepping is also applicable in the BSPF method.
\end{highlights}


\end{frontmatter}



\section{Introduction}
\label{sec1}

The spectral method is widely applied to solve partial differential equations (PDEs) with high fidelity across different fields from oceanography to plasma physics \citep{ducrozet2016hos,goerler2011global}. Orszag \cite{orszag1971numerical} first introduced the classical Fourier pseudo-spectral method for viscous flow simulations, leveraging the Fourier transform to achieve spectral accuracy. With modern high-performance Fast Fourier Transform (FFT) implementations \citep{frigo1998fftw}, the Fourier method delivers both high numerical accuracy and computational efficiency, making it a popular choice for PDE problems with simple geometries and periodic boundaries \cite{boyd2001chebyshev}.

For non-periodic problems, the Fourier method suffers from the Gibbs phenomenon at the boundaries, causing numerical instabilities or divergence \cite{gottlieb1997gibbs}. Gottlieb and Orszag \cite{gottlieb1977numerical} addressed this issue by using non-trigonometric basis functions such as Chebyshev, Legendre, or Jacobi polynomials, which are orthogonal, complete, and better suited for non-periodic problems. The Chebyshev pseudo-spectral method is particularly efficient for problems requiring high resolution at the boundary, thanks to its edge-clustering nodes \cite{GibsonHalcrowCvitanovicJFM08,riva2017uncertainty,boyd2014dynamics}. It can also leverage efficient FFT libraries through the application of the discrete cosine transform (DCT) \cite{gottlieb1977numerical} in its implementation. As a result, the Chebyshev method forms the cornerstone for multi-domain spectral methods \cite{pinelli1994chebyshev} and spectral element methods \cite{karniadakis2005spectral} with greater flexibility.

Despite its success, the reliance of the Chebyshev method on the edge-clustering grids also has its own disadvantages. For problems where critical physical phenomena are away from boundaries, such as shock waves \cite{don1994numerical}, shear instabilities \cite{chen2011accurate}, and internal waves \cite{hibiya1996direct}, it is less efficient to solve them with an edge-clustering grid. Techniques like domain decomposition \cite{pinelli1994chebyshev} or local refinement \cite{hashemi2017chebfun} might mitigate this issue, but they require \textit{a priori} assumptions on the solution and also introduce extra computational costs. Furthermore, the small grid of the Chebyshev points near the boundary imposes a stiff CFL condition that requires very small time steps and computationally intensive implicit time schemes  \cite{boyd2001chebyshev}.

An alternative approach to achieve spectral solutions for non-periodic problems is to extend non-periodic functions into periodic ones and keep using the Fourier-based methods. Early work by Gazdag \cite{gazdag1973numerical} used a cosine transition to smoothly connect domain endpoints, allowing Fourier methods to solve the viscous Burgers' equation. This Fourier extension concept has since been expanded with continuation functions based on Hermite \cite{boyd2002comparison}, Taylor \cite{fu2012modified}, Gram polynomials \cite{bruno2010high, albin2011spectral, bruno2022two}, and radial basis functions \cite{fryklund2018partition}. However, these methods often require the expansion of the domain and increase the computational cost, particularly for high-dimensional problems \cite{fu2012modified}.


Researchers have also explored the possibility of using an auxiliary function that matches the boundaries of the original function to achieve periodicity directly within the original domain. Roache \cite{roache1978pseudo} proposed a Taylor polynomial to enforce endpoint continuity of the original function and its derivatives. However, high-order polynomial matching (up to the 7th order \cite{roache1978pseudo}) is well-known to be susceptible to the Runge’s phenomenon, causing large oscillations and errors at the boundaries \cite{trefethen2019approximation}. In the subsequent studies, researchers explored the application of alternative functions, e.g., Gegenbauer polynomials \cite{vozovoi1996analysis} and Heaviside functions \cite{rim2013gibbs}. While these approaches offer better regularity than Taylor polynomials, they are either global approximants (still prone to extra oscillations) or low in smoothness, both of which hinder the spectral accuracy of the original Fourier method.

In parallel to spectral method developments, the field of computational mechanics has witnessed significant advances through Isogeometric Analysis (IGA), introduced by Hughes et al. \cite{hughes2005isogeometric}. The main novelty of IGA is the application of the B-spline basis functions employed in both geometry representation and finite-element analysis \cite{cottrell2009isogeometric}. The key advantage of B-splines lies in their inherent smoothness properties: B-splines of degree $p$ provide $C^{p-1}$ continuity \cite{de1978practical}, offering superior regularity compared to traditional $C^0$ finite element basis functions. This high-order continuity, combined with their local support and optimal approximation properties, enables IGA to be applied across diverse applications from structural mechanics to fluid dynamics \cite{cottrell2006isogeometric,tagliabue2014isogeometric, gondegaon2016static}.

Drawing inspiration from these superior properties of B-splines in IGA, we propose a novel B-spline-periodized Fourier (BSPF) pseudo-spectral method for non-periodic problems. The BSPF method provides a flexible framework to achieve and regularize the periodization of non-periodic functions with high-order accuracy and small computational overhead. The B-spline part is analytically evaluated, and the periodic residual is treated with the Fourier method. As a result, it can be performed on uniform grids or analytically mapped to grids with specific features while keeping the near-spectral convergence, making it well-suited for problems where critical physics are away from the boundaries and are inefficient to be resolved by the Chebyshev grid.

The remainder of this paper is structured as follows. Section 2 introduces the BSPF formulation and analyzes its theoretical errors and computational complexity. Section 3 explores the numerical behavior of the differentiation and integration operations on analytical test functions. Section 4 demonstrates its application to two nonlinear PDE problems: simulating the one-dimensional Burgers' equation and two-dimensional shallow-water equations. Section 5 concludes with a summary of the main properties and advantages of the BSPF method and discusses possible future extensions.

\section{B-spline-periodized Fourier method}
\subsection{B-spline periodization}
Given a non-periodic function $f(x)$ defined on a uniform grid $x_j$ within a closed interval $[a,b]$:
\begin{equation}
x_j = a + j \Delta x, \quad \Delta x = \frac{b-a}{N - 1}, \quad j = 0, 1, \dots, N-1, 
\label{eq:grid}
\end{equation}
Our objective is to compute its derivative $f'(x)$ by combining a B-spline approximation with a Fourier-based correction. The original function $f(x)$ is decomposed to
\begin{equation}
f(x_j) = f_s(x_j) + r(x_j),  
\label{eq:fx}
\end{equation}
where $f_s(x_j)$ is a non-periodic B-spline approximation that satisfies boundary constraints up to certain orders of derivatives, and $r(x_j)$ is a periodic residual which can be treated with the classical Fourier spectral method.

\subsubsection{B-spline approximation}
We approximate $f_s(x_j)$ in Eq.  \ref{eq:fx} using a set of B-spline basis functions $B_{i,p}(x)$ of degree $p$ \cite{de1978practical}:
\begin{equation}
f_s(x_j) = \sum_{i=0}^{n-1} P_i B_{i,p}(x_j),
\label{eq:fsx}
\end{equation}
B-splines are piecewise polynomial functions renowned for their smoothness and local support, making them well-suited for stable and efficient function approximation. where $P_i$ are the control coefficients to be determined; \( n \) is the number of basis functions determined by a knot vector $\mathbf{z}$ distributed within $[a,b]$:
\begin{equation}
\mathbf{z} = \{ \zeta_0, \zeta_1, \dots, \zeta_{n+p} \}.
\end{equation}
We first choose a knot distribution with clamped (repeated) endpoints to enforce boundary conditions within $[0,1]$:
\begin{equation}    
\mathbf{z} = \left\{
\underbrace{0, 0, \dots, 0}_{p+1},
\zeta\left(\frac{1}{n-p}\right), \dots,\zeta\left(\frac{n-p-1}{n-p}\right),
\underbrace{1, 1, \dots, 1}_{p+1}
\right\},
\label{eq:Xi}
\end{equation}
and then scale it to \([a, b]\). The repeated knots at the boundaries, each with  \(p+1\) multiplicity, ensure that the spline interpolates the endpoints \cite{de1978practical}.  The interior knot distribution can be a uniform one ($\zeta(x) = x$), and it can also be further mapped to an edge-clustering one to provide extra regularization of boundary matching. In this work, we choose a $\tanh$ distribution with a control parameter $\beta$ for the edge-clustering mapping
\begin{equation}
    \zeta(x) = \frac{\tanh \beta x}{\tanh \beta}.
    \label{eq:knots_dist}
\end{equation}
Note that Eq. \ref{eq:knots_dist} is just one example of edge-clustering distributions, and one can freely choose other distributions, such as a polynomial or a Chebyshev-like cosine function. We leave the question of the optimal knot distribution to future explorations and choose Eq. \ref{eq:knots_dist} here due to its simplicity. 

The basis function $B_{i, p}(x)$ can be computed recursively using the Cox-de Boor formula \cite{de1978practical}:
\begin{equation}
B_{i, p}(x)=\frac{x-\zeta_i}{\zeta_{i+p}-\zeta_i} B_{i, p-1}(x)+\frac{\zeta_{i+p+1}-x}{\zeta_{i+p+1}-\zeta_{i+1}} B_{i+1, p-1}(x).
\label{eq:Bkm}
\end{equation}
When $B_{i, p}(x)$ is recursively reduced to $B_{i, 0}(x)$, it becomes a rectangular window function:
\begin{equation}
    B_{i, 0}(x)= \begin{cases}1 & \zeta_i \leq x< \zeta_{i+1}, \\ 0 & \text { otherwise. }\end{cases}
\label{eq:Bk0}
\end{equation}

\subsubsection{Boundary constraint}
To achieve a smooth periodization of $f(x)$ with B-spline of degree $p$, we enforce the following boundary constraint up to the $(p-1)$th derivatives at the endpoints:
\begin{equation}
f_s^{(k)}(a) = f^{(k)}(a), \quad f_s^{(k)}(b) = f^{(k)}(b), \quad k = 0, 1, \dots, p-1.
\label{eq:boundary_matching}
\end{equation}
Substituting Eq. \ref{eq:boundary_matching} into Eq. \ref{eq:fsx}, we construct the following linear system:
\begin{equation}
C_{k,i}P_i = d_k,
\label{eq:CPd}
\end{equation}
where $C_{k,i} \in \mathbb{R}^{2p \times n}$ is the boundary constraint matrix:
\begin{equation}    
C_{k,i} = 
\begin{cases}
B_{i,p}^{(k)}(a), & 0 \le k < p, \\
B_{i,p}^{(k-p )}(b), & p \le k < 2p,
\end{cases}
\end{equation}
and
\begin{equation} 
d_k = [f^{(0)}(a), \dots, f^{(p-1)}(a), f^{(0)}(b), \dots, f^{(p-1)}(b)]^T.
\end{equation}
where $d_k \in \mathbb{R}^{2p}$ is the boundary derivative vector to compute from the data. 

The linear system defined in Eq. \ref{eq:CPd} is the centerpiece of the BSPF method. Since the boundary constraint matrix $C$ only contains the information of B-spline basis functions, it can be precomputed only once without the information of $f(x)$. Eq. \ref{eq:CPd} also offers a natural way to accommodate boundary conditions. For a Dirichlet condition, we simply replace the selected rows of the boundary-derivative vector $d_k$ with prescribed values, e.g.\ $d_0=f(a)$ and $d_{p}=f(b)$. A Neumann condition is obtained by modifying the first-derivative rows, e.g.\ $d_1=f'(a)$ and $d_{p+1}=f'(b)$.

\subsubsection{Boundary derivatives}
The boundary derivatives $d^L_k=f^{(k)}(a)$ and $d^R_k = f^{(k)}(b)$ appearing in the RHS of Eq. \ref{eq:CPd} are solved from a $m$-point ($m\geq p$) one-side stencil $f(x_j)$ using the Taylor expansion:
\begin{equation}
y^L_j= f(x_j) \approx \sum_{k=0}^{m-1} \frac{(x_j - a)^k}{k!} f^{(k)}(a), \quad j=0,\dots,m-1,
\end{equation}
\begin{equation}
y^R_j = f(x_{N-1-j}) \approx \sum_{k=0}^{m-1} \frac{(x_{N-1-j} - b)^k}{k!} f^{(k)}(b), \quad j=0,\dots,m-1.
\end{equation}
In matrix form, it can be expressed as:
\begin{equation}    
 A^{L}_{j,k} d^{L}_k = y^{L}_{j}, \quad  A^{R}_{j,k} d^{R}_k = y^{R}_{j} 
\end{equation}
where $A_{j, k}^{L}, A_{j,k}^{R} \in \mathbb{R}^{m \times m}$ are the Vandermonde matrices corresponding to the left and right boundaries:
\begin{equation}
    A_{j, k}^{L} = \frac{(x_j - a)^k}{k!}, \quad A_{j,k}^{R} = \frac{(x_{N-1-j} - b)^k}{k!}
\end{equation}
In the case of a uniform grid, they can be reduced to
\begin{equation}
    A_{j, k}^{L} = -A_{j,k}^{R}= \frac{(j\Delta x)^k}{k!}.
\end{equation}

\subsubsection{Regularization}
To make the linear system defined in Eq. \ref{eq:CPd} solvable, we need to at least have the number of basis functions $n\geq 2 p$. The system is determined when $n = 2p$ and under-determined when $n > 2p$. When the original function has large gradients at the boundaries, a determined system could be ill-conditioned, leading to a $f_s$ with huge oscillations. Therefore, it is usually necessary to introduce extra basis functions to regularize $f_s$.

Since $C_{k,i} P_i = d_k$ is under-determined when $n > 2p$, we select the coefficients $P_i$ that minimize the $L^2$ error to $f(x)$:
\begin{equation}
\min_{\mathbf{P}} \int_a^b \left( \sum_{i=0}^{n-1} P_i B_{i,p}(x) - f(x) \right)^2 dx + \lambda \|P_i\|^2,    
\end{equation}
while still subject to $C_{k,i} P_i = d_k$. $\lambda \ge 0$ is the Tikhonov regularization parameter.

Let $D \in \mathbb{R}^{n \times N}$ be the matrix of basis function evaluations, $D_{i,j} = B_{i,p}(x_j)$, and $W$ be the diagonal matrix of trapezoidal quadrature weights
\begin{equation}
W = \mathrm{diag}(w_0, w_1, \dots, w_{N-1}),
\end{equation}
where
\begin{equation}
w_j =
\begin{cases}
\frac{h}{2}, & j = 0 \ \text{or}\ j = N-1, \\
h, & 1 \le j \le N-2.
\end{cases}
\end{equation}
This construction ensures that the discrete $L^2$ inner product
\begin{equation}
\langle u, v \rangle_{L^2} \approx u^T W v
\end{equation}
approximates the continuous integral. We further define
\begin{equation}
Q = D W D^T, \quad \mathbf{L} = D W \mathbf{f}.
\end{equation}
The resulting Karush-Kuhn-Tucker (KKT) system is
\begin{equation}
\begin{bmatrix}
2(Q + \lambda I) & -C^T \\
C & 0
\end{bmatrix}
\begin{bmatrix}
\mathbf{P} \\
\boldsymbol{\mu}
\end{bmatrix}
=
\begin{bmatrix}
2 \mathbf{L} \\
\mathbf{d}
\end{bmatrix},
\end{equation}
where $\boldsymbol{\mu}$ are Lagrange multipliers, $\mathbf{P}$ is the vector of the regularized spline coefficients $P_i$, $\mathbf{d}$ is the vector of the boundary derivatives $d_k$, and $\mathbf{f}$ is the function sample vector $f(x_j)$. Solving this system yields the regularized spline coefficients $\mathbf{P}$, which defines a B-spline function $f_s(x)$ that
\begin{enumerate}
    \item matches the boundary values and derivatives of $f(x)$ up to $p$-th order;
    \item is the closest approximation to the original function $f(x)$ in the $L_2$ sense in the function space 
\(\mathcal{V}_n\)
spanned by the $n$ B-spline basis functions.
\end{enumerate}

\subsection{Basic operations}
\subsubsection{Differentiation}
Once we have obtained the B-spline periodization of the original function $f(x)$, we can take advantage of the property that B-spline basis functions $B_{i,p}(x)$ can be analytically differentiated to compute the derivatives of $f(x)$. The B-spline derivative is given by
\begin{equation}
f_s'(x) = \sum_{i=0}^{n-1} P_i B_{i,p}^{(1)}(x).
\end{equation}
which $B_{i,p}^{(1)}(x)$ can be recursively evaluated by \cite{de1978practical}: 
\begin{equation}
 B^{(1)}_{i,p}(x) = \frac{p}{\zeta_{i+p} - \zeta_i} B_{i,p-1}(x) - \frac{p}{\zeta_{i+p+1} - \zeta_{i+1}} B_{i+1,p-1}(x).
\end{equation}
The periodic correction term $r'(x)$ is obtained by the classical Fourier spectral differentiation:
\begin{equation}
r'(x_j) = \mathcal{F}^{-1} \left[ i \, \omega \, \mathcal{F}(f - f_s) \right],
\end{equation}
where $i$ is an imaginary unit, $\omega$ is the frequency vector of the FFT. The final derivative is
\begin{equation}
f'(x) = f_s'(x) + r'(x).
\end{equation}
\subsubsection{Antiderivative (indefinite integration)}
The BSPF method can also be applied to compute the antiderivative, or indefinite integral of $f(x)$. We still obtain a boundary-matching B-spline $f_s(x)$. The antiderivative of $f_s(x)$ can also be evaluated by 
\begin{equation}
\mathcal{I}[f_s(x)] = \sum_{i=0}^{n-1} P_i \int B_{i,p}(x),
\end{equation}
where the antiderivative of $B_{i,p}(x_j)$ is computed analytically by \cite{schumaker2007spline}:
\begin{equation}
\int B_{i, p}(x_j) d x= \frac{\zeta_{i+p+1}-\zeta_i}{p+1} B_{i, p+1}(x)+S_i
\end{equation}
The periodic correction term $r(x)$ is computed using the standard Fourier spectral approach:
\begin{equation}
\mathcal{I}[r(x_j)] = \mathcal{F}^{-1} \left[ \frac{1}{i \, \omega} \, \mathcal{F}(f - f_s) \right],
\end{equation}
where $\omega$ is the frequency vector of the FFT. The final integrated function is
\begin{equation}
\mathcal{I}[f(x)] = \mathcal{I}[f_s(x)] + \mathcal{I}[r(x)] + S.
\end{equation}
where the integration constant $S$ can be determined by boundary conditions.

\subsubsection{Grid mapping \label{sec:mapping}}
The BSPF method can also be easily applied to non-uniform grids through analytical mapping without losing its accuracy. Let $x_j$ be the uniform grid points defined on the interval $[a,b]$. We introduce $\xi_j = g(x_j)$ as the mapped non-uniform grids, where $g(x)$ is a function whose derivatives $g'(x)$ can be analytically computed, and $g(x)$ also maps the endpoints onto themselves ($g(a) = a, g(b) = b$). We further define the compound function $F(x) = f(g(x))$ on the interval $[a,b]$. Based on the chain rule of differentiation, we have the following relation:
\begin{equation}
    \frac{df}{d\xi} = \frac{F'(x)}{g'(x)},
\end{equation}
As we mentioned before, $g'(x)$ can be analytically evaluated, and $F'(x)$ can be computed from the BSPF method on the uniform grid $x_j$. As long as $g'(x_j) \neq 0$ within $[a,b]$, we can easily retrieve the derivatives of the original function $f(x)$ on the non-uniform grid $\xi_i$. Similarly, the integration on the mapped grid is
\begin{equation}
    \int{f(x)}{d\xi} = \int{F(x)}{g'(x)}dx +S,
\end{equation}
where the integration constant $S$ can be determined by boundary conditions.
\subsection{Algorithmic complexity\label{sec:time}}
The key steps of the BSPF method and their time complexities are summarized in Algorithm \ref{algo:bfpsm}. For the sake of conciseness, here we demonstrate the differentiation operation with a determined system ($n = 2p$) with a $p$-point boundary stencil ($m=p$), which makes the boundary constraint matrix $C$ a square matrix ($C \in \mathbb{R}^{2p\times 2p}$):

\begin{algorithm}[h!]
\caption{BSPF (Differentiation)}
\label{algo:bfpsm}
\begin{algorithmic}
\Require Grid $\mathbf{x}\in\mathbb R^N$; function values $\mathbf f\in\mathbb R^N$; B-spline degree $p$; knots $\mathbf{z}\in\mathbb R^{3p+1}$; Grid number $N$;

\State \textbf{Pre-compute (using $\mathbf{x}$ and $\mathbf{z}$):}
\State Build a basis function $B,B'$ on $\mathbf{x}$. \Comment{$\mathcal O(pN)$}
\State Assemble a boundary constraint matrix $C$. \Comment{$\mathcal O(p^2)$}
\State Form and invert Vandermonde matrices $A^{L},A^{R}$. \Comment{$\mathcal O(p^3)$}
\State LU factorization of $C$. \Comment{$\mathcal O(p^3)$}

\State \textbf{Apply (using $\mathbf f$):}
\State Compute boundary targets $\mathbf d$. \Comment{$\mathcal O(p^2)$}
\State Solve pre-factorized $C\,\mathbf P=\mathbf d$. \Comment{ $\mathcal O(p^2)$}
\State Compute $\mathbf f_s=B^\top\mathbf P$ and $\mathbf f'_s={B'}^\top\mathbf P$. \Comment{$\mathcal O(pN)$}
\State Periodic correction: $\widehat{\mathbf r'}=\mathrm{irfft}(i\,\boldsymbol\omega\odot \mathrm{rfft}(\mathbf f-\mathbf f_s))$. \Comment{$\mathcal O(N\log N)$}
\State \textbf{Return} $\widehat{\mathbf f'}=\mathbf f'_s+\widehat{\mathbf r'}$.
\end{algorithmic}
\end{algorithm}

Due to the constraint of the Taylor derivative approximation at the boundary, the practical choice of the degree of the basis function $p$ is $p\sim 10$. Therefore, we can assume that $N \gg p$.  As a result, the leading term of the time complexity in the precomputing stage is the construction of basis functions $B,B'$ ($\mathcal O(pN)$) and the inversion of $A^{L},A^{R}$ ($\mathcal O(p^3)$) . They only need to be computed once when the grid and knot vector are fixed.

In the apply stage, in which the online information of $\mathbf{f}$ is used, the leading terms of time complexity are the B-spline derivatives ($\mathcal O(pN)$) and the FFT-based periodic correction ($\mathcal O(N\log N)$). With the increase of $N$, the time complexity is dominated by the FFT-related operation. Therefore, the overall time complexity of the BSPF algorithm is asymptotically at the same level as the classical Fourier or Chebyshev pseudo-spectral methods. The integration operation shares the same time complexity as the differentiation, as it is also dominated by FFT.
\subsection{Theoretical error analysis}\label{subsec:theoretical-error}
In this section, we further analyze the theoretical error of the BSPF method on a uniform grid with \(N\) points over \([a,b]\). Similar to the complexity analysis, we consider a determined system $n = 2p$. We approach this problem by a two-step approach, assuming we have the exact boundary derivatives first and then injecting the errors induced by the finite difference approximation. 

If the endpoint derivatives up to order \(p-1\) are exact, the residual \(r^{\mathrm{ex}}\) closes periodically with \(C^{p-1}\) regularity at the boundaries. According to the property of Fourier differentiation \cite{trefethen2000spectral}, after applying \(\nu\)-th pseudo–spectral differentiation, the discarded tail obeys
\begin{equation}
\big\|D_{N}^{(\nu)} r^{\mathrm{ex}} - (r^{\mathrm{ex}})^{(\nu)}\big\|
\;=\;\begin{cases}
\mathcal O\!\big(N^{-(\,p-\nu\,)}\big), &( C^{p-1}),\\[0.25em]
\mathcal O\!\big(N^{\nu} e^{-\alpha N}\big), & \text{(analytic)}.
\end{cases}
\end{equation}
In practice, the endpoint derivatives are computed from an \(m\)-point one–sided FD stencil. For the \(k\)-th derivative (\(k=0,\dots,p\)) this gives
\begin{equation}
\delta d_k=\widehat d_k-d_k=\mathcal O(\Delta x^{\,m-k})=\mathcal O\!\big(N^{-(m-k)}\big).
\end{equation}
Because the periodization matches \(\widehat d_k\), the actually used residual \(r\) inherits boundary defects
\begin{equation}
J_k:=r^{(k)}(a)-r^{(k)}(b)=-\delta d_k=\mathcal O\!\big(N^{-(m-k)}\big).
\end{equation}
For the component carrying a jump in the \(k\)-th derivative one has \(|\widehat r_{1,\ell}|\le C\,|J_k|\,|\ell|^{-(k+1)}\). After \(\nu\)-th differentiation and truncation,
\begin{equation}
\sum_{|\ell|>N} |\ell|^{\nu}\,|\widehat r_{1,\ell}|
\;\lesssim\; |J_k|\,\sum_{|\ell|>N} |\ell|^{-(k+1-\nu)}
\;=\; \mathcal O\!\big(|J_k|\,N^{-(\,k-\nu\,)}\big),
\end{equation}
so with \(|J_k|=\mathcal O(N^{-(m-k)})\) each channel contributes \(\mathcal O(N^{-(m-\nu)})\) . Therefore the FD–closure floor for the \(\nu\)-th derivative is
\begin{equation}
\big\|D_{N}^{(\nu)} r - r^{(\nu)}\big\| \;=\; \mathcal O\!\big(N^{-(\,m-\nu\,)}\big).
\end{equation}
Combining the two and adding an \(N\)-independent noise term \(\sigma\) from rounding errors, we obtain the differentiation error bound
\begin{equation}
\big\|D_{N}^{(\nu)} f - f^{(\nu)}\big\|
\;\lesssim\;
\underbrace{\mathcal O(N^{\nu} e^{-\alpha N})}_{\text{analytical tail }}
\;+\;
\underbrace{\mathcal O(N^{-(\,p-\nu\,)})}_{\text{algebraic tail}}
\;+\;
\underbrace{\mathcal O(N^{-(\,m-\nu\,)})}_{\text{Taylor approx.}}
\;+\;
\underbrace{\sigma}_{\text{round-off err.}}\,.
\end{equation}
In some cases, the error introduced by the finite difference approximation at the boundary is greater than the Fourier truncation, so we might only observe the $m-\nu$ algebraic convergence before it reaches the round-off error floor.

Similarly, for the integration operation we have
\begin{equation}
\big\|I_{N}^{(\nu)} f - f^{(\nu)}\big\|
\;\lesssim\;
\underbrace{\mathcal O(N^{-\nu} e^{-\alpha N})}_{\text{analytical tail }}
\;+\;
\underbrace{\mathcal O(N^{-(\,p+\nu\,)})}_{\text{algebraic tail}}
\;+\;
\underbrace{\mathcal O(N^{-(\,m+\nu\,)})}_{\text{FD approx.}}
\;+\;
\underbrace{\sigma}_{\text{round-off err.}}\,.
\end{equation}

\section{Numerical results: differentiation and integration}
\subsection{Differentiation\label{sec:diff}}
\begin{figure}[h!]
\centering
\includegraphics[width=1\linewidth]{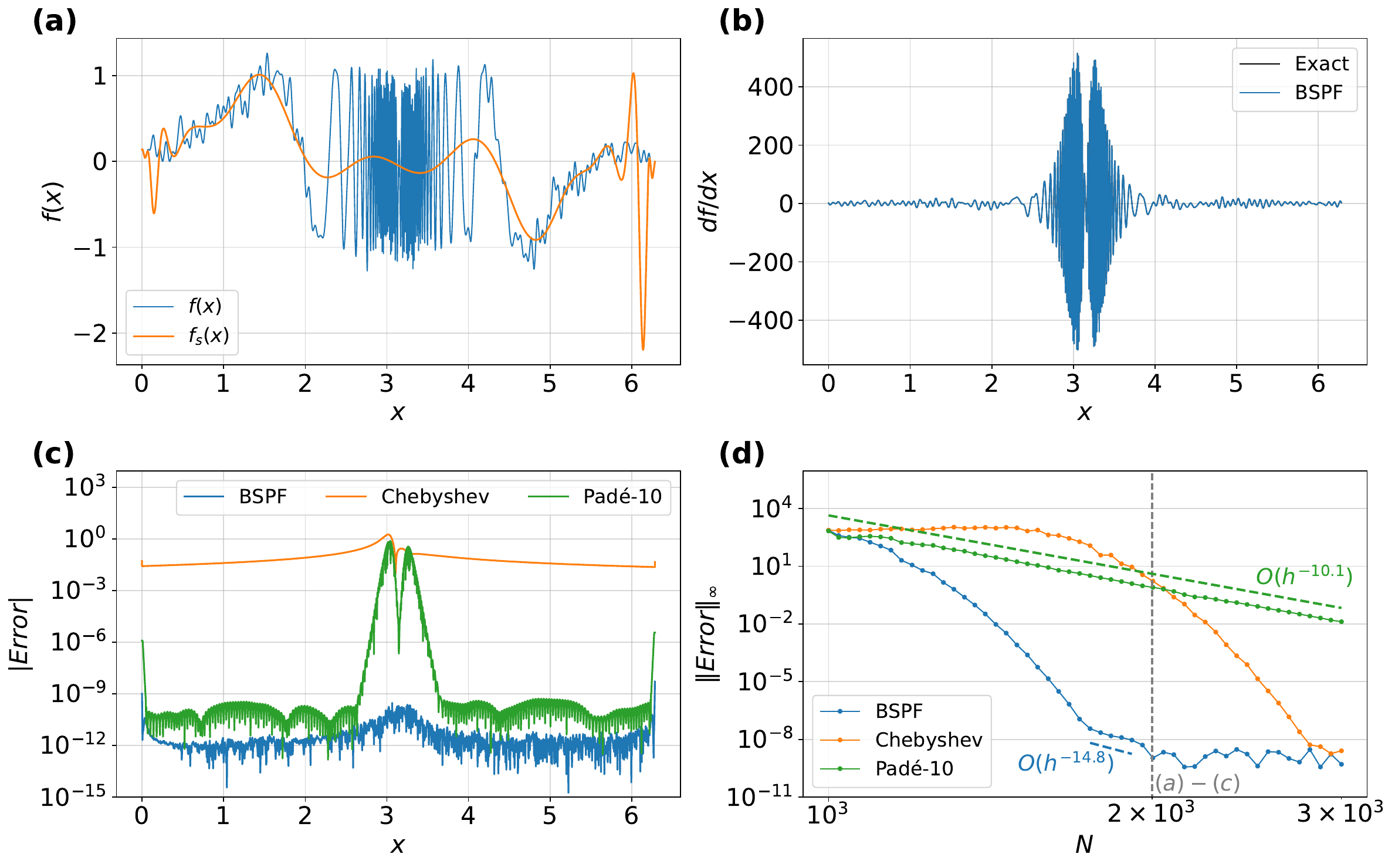}
\caption{\label{fig:1D-diff}Computing the first-order derivative of a locally oscillatory function with noise: (a) the original function $f(x)$ (blue) and its non-periodic B-spline component $f_s(x)$ (orange); (b) the exact derivative (black) and the BSPF numerical result (blue); (c) error distributions to the exact $f'(x)$ by the BSPF (blue), Chebyshev (orange), and the 10th-order compact finite difference (green) methods; (d) grid convergences of the $L^\infty$ norm error of the aforementioned methods. The results of (a)-(c) are computed on the grid $N = 2000$. }
\end{figure}

We consider a locally oscillatory function within the interval $[0,2\pi]$ with a non-periodic random noise (Figure \ref{fig:1D-diff}a) to examine the numerical performance of the BSPF method in the task of derivative computing: 
\begin{equation}
    f(x) = \sin \left(\frac{x}{1.02+\cos(x)}\right) + \epsilon(x).
    \label{eq:fluc}
\end{equation}
The synthetic noise signal $\epsilon(x)$ is generated by randomly sampling its frequencies $\kappa_i \sim \mathcal{U}(1,100)$, amplitudes $\mathcal{A}_i\sim \mathcal{U}(0,0.01)$ and phases $\varphi_i\sim\mathcal{U}(0,2\pi)$, where $\mathcal{U}$ is the uniform distribution. We then construct $\epsilon(x)$ by: 
\begin{equation}
\epsilon(x)=\sum_{i=1}^{1000} \mathcal{A}_i \cos\!\big(\kappa_i\,x+\varphi_i\big).
\end{equation}
Since $\epsilon(x)$ is a linear superposition of cosine functions with prescribed parameters, it is still analytically differentiable. The BSPF method is applied to compute its first-order derivative $f'(x)$, and its results are compared to the results of the Chebyshev method, the 10th-order Padé-type compact finite difference scheme \cite{lele1992compact}, and the exact solution. The parameters of the BSPF method used in the numerical experiments are summarized in Table \ref{tab:parameters}.
\begin{table}[ht]
\caption{BSPF parameter values used in the numerical experiments.}
\label{tab:parameters}
\centering
\begin{tabular}{lcccc}
\toprule
Parameter & Degree ($p$) & Basis ($n$) & Border stencil ($m$)& Clustering ($\beta$) \\
\midrule
Values & 11& 44 & 16& 3.0 \\
\bottomrule
\end{tabular}

\end{table}

We choose an 11-degree B-spline ($p=11$) to match the boundary up to 10th derivatives. Due to the presence of random noise, we apply a regularized BSPF system with extra basis functions ($n = 4p$ instead of the minimal $2p$)  and an edge-clustered $\tanh$ knot distribution (Eq. \ref{eq:knots_dist} with $\beta=3.0$) to stabilize the B-spline fitting. The boundary derivatives are computed with a 16-point one-sided finite difference stencil ($m = 16$).  The BSPF and compact finite difference computations are carried out on a uniform grid with $N=2000$, and the Chebyshev computation is performed on a Chebyshev grid with the same grid number. To investigate the convergence of each method, we compute its $L^\infty$ norm errors of the numerical results to the analytical solution using grids ranging from $N=1000$ to $N=3000$.
 
Using these configurations, a B-spline fit $f_s(x)$ is constructed using Algorithm \ref{algo:bfpsm} to match the boundary derivatives of the original function $f(x)$ up to its 10th-order derivatives  (Figure \ref{fig:1D-diff}a), allowing the application of Fourier differentiation of the periodic residual $r(x)$. $f(x)$ is highly oscillatory around $x = 1.5\pi$, causing steep fluctuations in $f'(x)$ (Figure \ref{fig:1D-diff}b) and posing significant challenges for numerical differentiation. In the case of $N=2000$ (Figure \ref{fig:1D-diff}c), the result of the 10th-order Padé scheme shows local error peaks ($\sim10^{0}$) around the oscillatory part of $f'(x)$, while the Chebyshev results contains large global errors $(10^{-2}$ to $10^0)$ across the entire domain.  The BSPF method, on the other hand, suppresses the error below $10^{-9}$ except for the right boundary point. In terms of error convergence (Figure \ref{fig:1D-diff}d), the Chebyshev and BSPF methods both initially demonstrate exponential convergence. However, the Chebyshev converges more slowly than the BSPF method due to the fact that the Chebyshev grid is stretched towards the domain boundary. As a result, it is less effective to capture the oscillation in the middle of the domain compared to the uniform grid used by the BSPF method. As the error is reduced to $\sim10^{-8}$, it starts to be dominated by the boundary derivatives approximation, thus the BSPF method starts to converge algebraically with a rate of $O(m-1)$ until reaching the rounding error floor, consistent with the theoretical error we present in Section \ref{subsec:theoretical-error}. 

\begin{figure}[h!]
\centering
\includegraphics[width=1\linewidth]{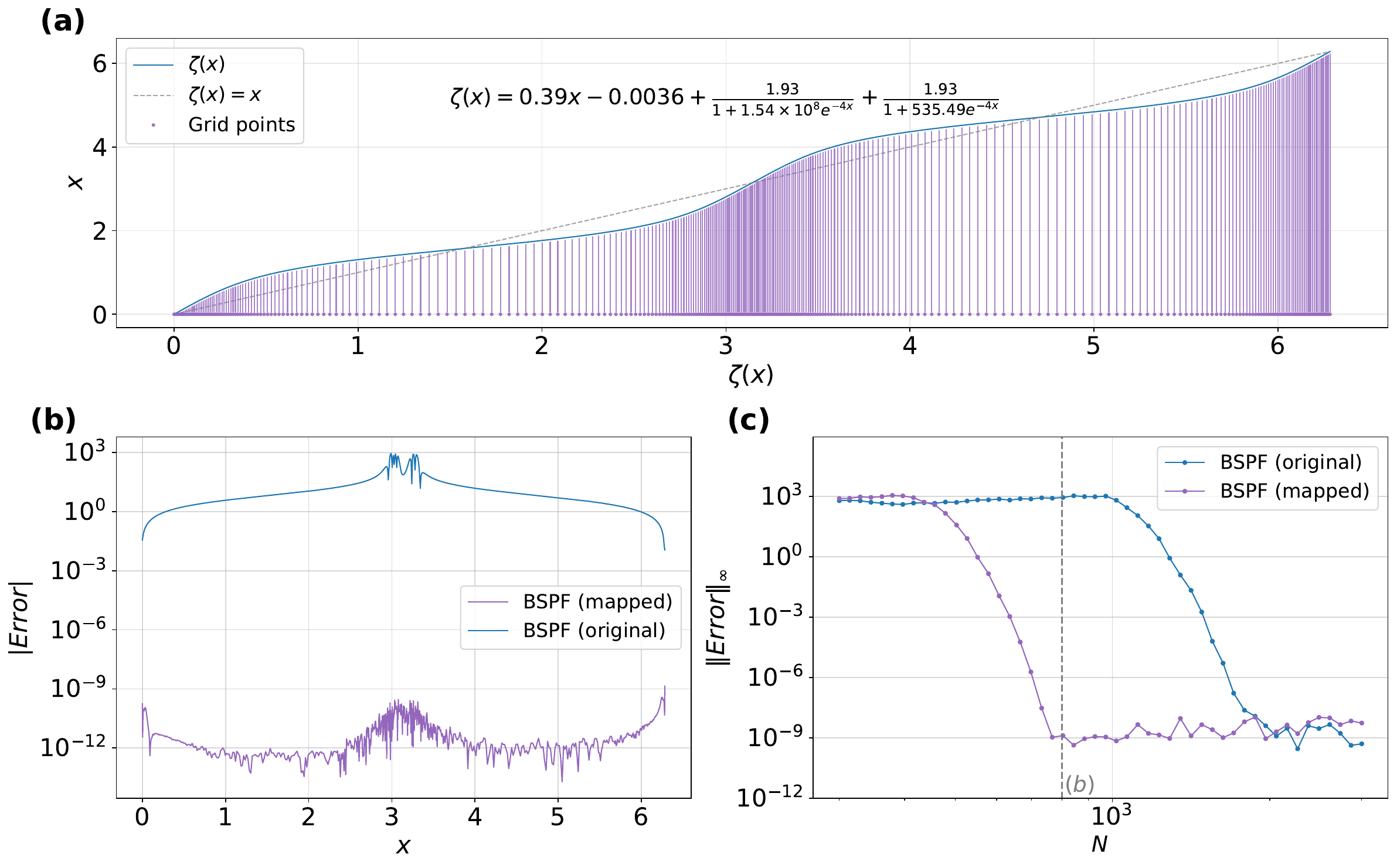}
\caption{\label{fig:1D-diff-refine}Computing the first-order derivative on a mapped grid: (a) analytical mapping function $\zeta(x)$ (blue) and mapped grid (purple); (b) Spatial distributions of the error magnitude of the BSPF method to the exact solution on the original uniform grid (blue) and mapped grid (purple); (c) grid convergences of the $L^\infty$ norm error of the BSPF on the uniform grid (blue) and on the mapped grid (purple). The result of (b) is computed on the grid $N = 800$.}
\end{figure}

As discussed in Section \ref{sec:mapping}, the analytical mapping can also be employed to apply the BSPF method on a non-uniform grid. Figure \ref{fig:1D-diff-refine}a illustrates an example of such a mapping function, $\zeta(x)$, constructed by smoothly concatenating two sigmoid functions to refine the mesh near the boundaries and at the center of the domain ($x=\pi$). On this mapped grid with local refinement ($N = 800$), the BSPF method achieves much lower errors ($< 10^{-9}$) than on a uniform mesh for computing the derivative of Eq. \ref{eq:fluc} (Figure \ref{fig:1D-diff-refine}b). Specifically, the error level of $\sim 10^{-9}$ is reached before $N=800$, owing to the enhanced resolution around $x=\pi$ and the improved boundary matching accuracy provided by the edge-clustered grid.

\begin{figure}[h!]
\centering
\includegraphics[width=0.62\linewidth]{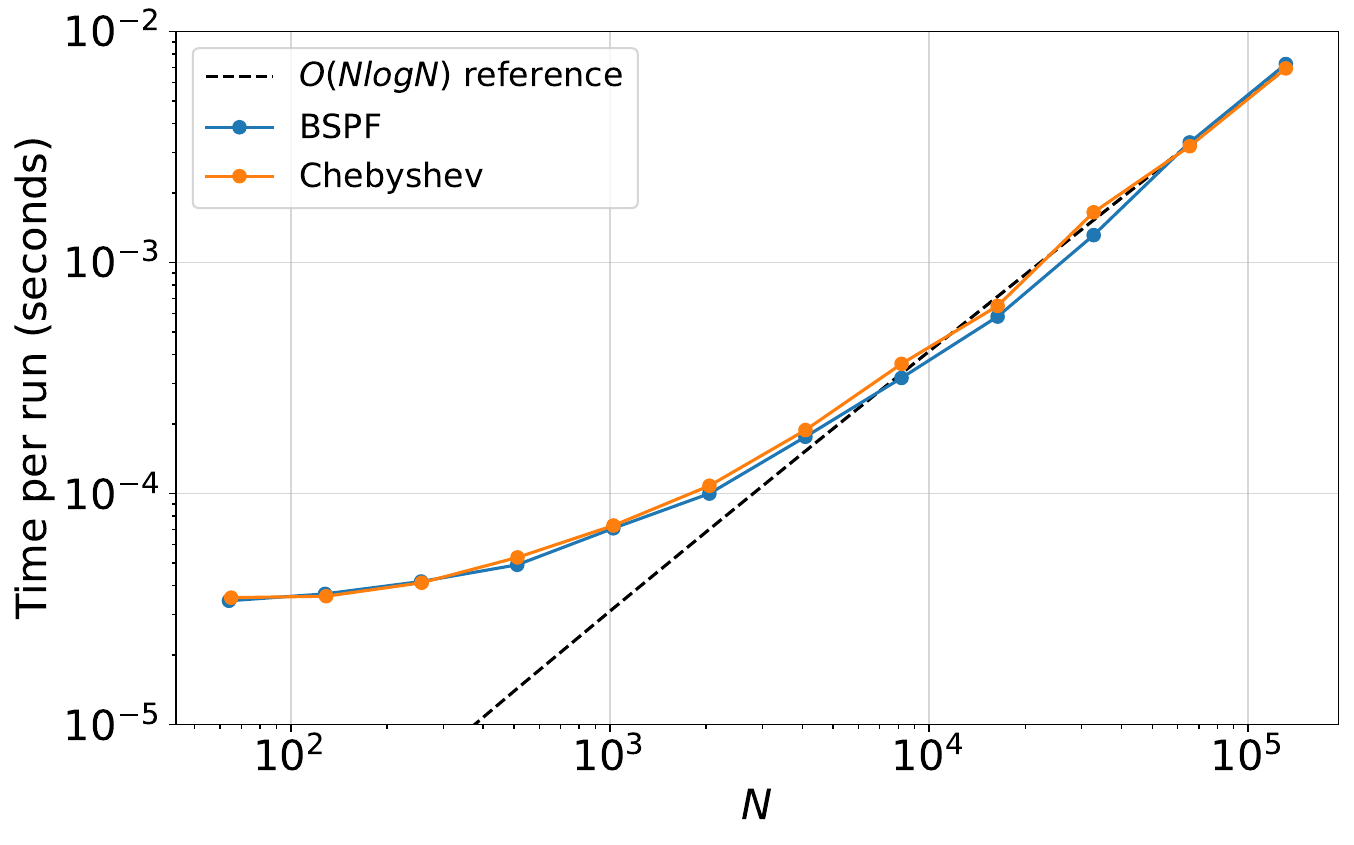}
\caption{\label{fig:time-complexity} Execution time of derivative computations using the BSPF (blue) and Chebyshev (orange) methods vs. the grid number $N$. The execution time is measured by averaging the wall time of 100 runs on an AMD EPYC 7453 server CPU.}
\end{figure}

Finally, Figure \ref{fig:time-complexity} compares the execution time of derivative computation using the BSPF and Chebyshev methods as a function of the problem size $N$. The results show that both methods exhibit nearly identical performance across a wide range of $N$, from 64 to 131072 grid points. Both curves closely follow the $\mathcal{O}(N \log N)$ reference trend line after the problem size approaches $10^4$, which is consistent with our time-complexity analysis in Section \ref{sec:time} showing that their computational complexity is dominated by the FFT-related operation for large problem sizes.

\subsection{Integration}
In the integration test, we consider a 1D ODE problem in the interval $[a,b]$ with Dirichlet boundary conditions: 
\begin{equation}
    \frac{d f(x)}{dx} = f'(x), \quad f(a) = u^{L}, f(b) = u^{R}. 
    \label{eq:poisson}
\end{equation}
In the 1D case, Eq. \ref{eq:poisson} can be solved by direct integration:
\begin{equation}
    f(x) = \int f'(x) dx + S,
\end{equation}
in which the boundary conditions are used to determine the integration constant $S$.

\begin{figure}[h!]
\centering
\includegraphics[width=1\linewidth]{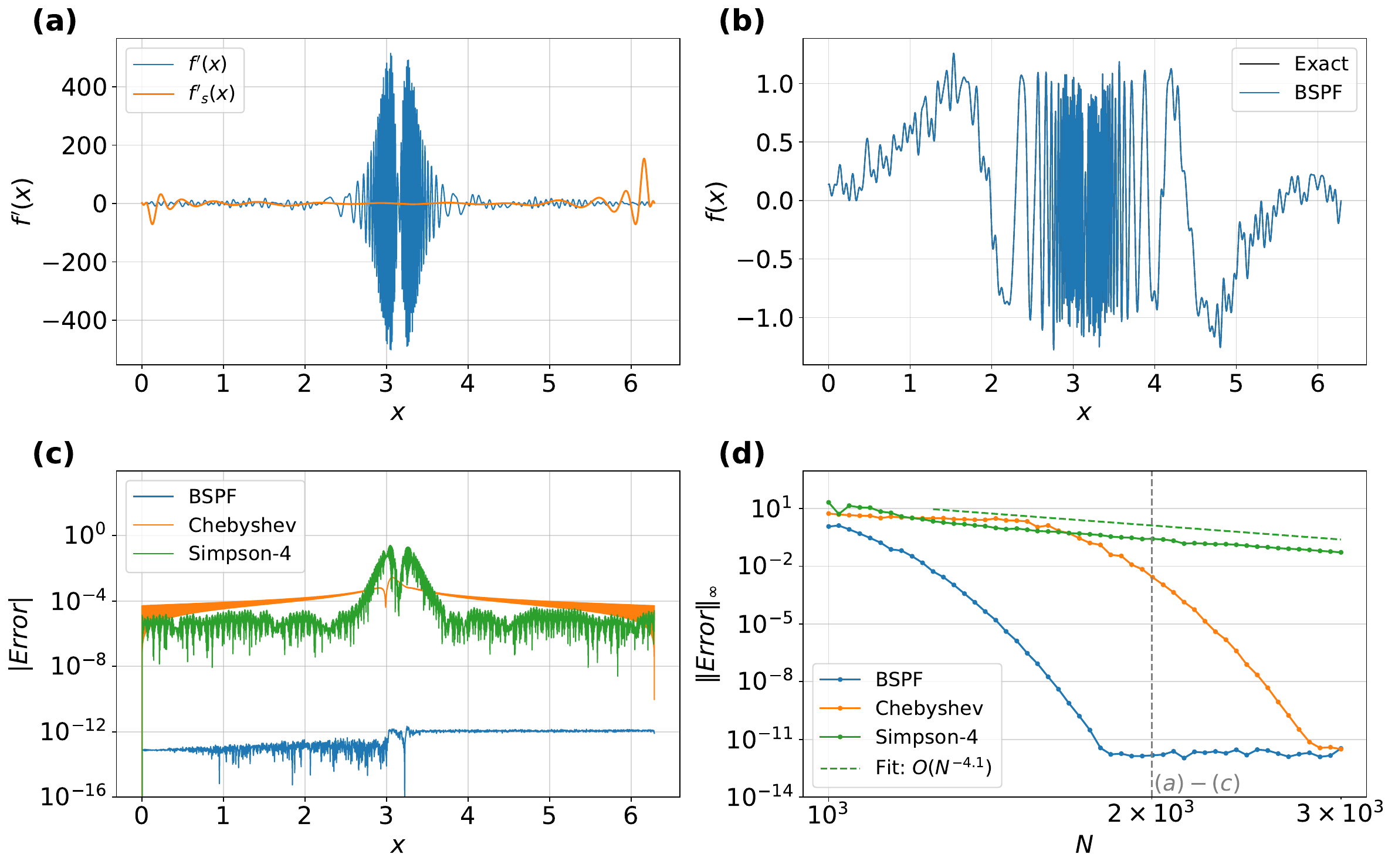}
\caption{\label{fig:poisson} Solving a 1D ODE problem through direct integration: (a) the original function $f'(x)$ (blue) and its non-periodic B-spline component $f'_s(x)$ (orange); (b) comparisons of the exact (black) and computed $f(x)$ (blue); (c) distribution of the absolute error to the exact solution; (d) grid convergence of the $L^2$ norm error of the BSPF (blue), Chebyshev  (orange) and the forth-order Simpson methods (green). }
\end{figure}

We consider the inverse problem of the differentiation discussed in Section \ref{sec:diff} and take the analytical derivative of Eq. \ref{eq:fluc} as the RHS of Eq. \ref{eq:poisson}. The numerical configuration is identical to the test case in Section \ref{sec:diff}. The results of integration are shown in Figure \ref{fig:poisson}.

As shown in Figure~\ref{fig:poisson}a, the input $f'(x)$ exhibits strong oscillations caused by local fluctuations in $f(x)$ (Figure~\ref{fig:poisson}b). Similar to the differentiation case, both the BSPF and Chebyshev methods achieve exponential convergence. However, due to the large variations of $f(x)$ in regions where the Chebyshev nodes are sparse, the Chebyshev method yields larger errors than its BSPF counterpart at the same grid resolution (Figure~\ref{fig:poisson}c) and converges significantly more slowly (Figure~\ref{fig:poisson}d). Meanwhile, the Simpson method converges at a much slower algebraic rate of $O(N^{-4})$.

\section{Numerical results: solving PDEs}
\subsection{Burgers' equation}
The Burgers' equation is an important non-linear PDE that arises in applied mathematics and fluid mechanics. It captures both nonlinear advection and viscous diffusion, making it a valuable test model for examining the numerical performance of the BSPF in the task of solving PDE.

The Burgers' equation in 1D can be expressed as:
\begin{equation}
    \frac{\partial u}{\partial t} + u \frac{\partial u}{\partial x} = \nu \frac{\partial^2 u}{\partial x^2},
    \label{eq:1d_burgers}
\end{equation}
where $u(x,t)$ is a time-dependent scalar function; $\nu$ is the viscosity. In certain initial conditions, Eq. \ref{eq:1d_burgers} can be solved exactly with the method of characteristics. The smooth step is one of those cases where the analytical solution of $u(x,t)$ is a traveling wave front:
\begin{equation}
    u_{exact}(x, t)=\frac{a+b+(b-a) e^{\eta}}{1+e^{\eta}}, 
    \quad   \eta=\frac{a}{\nu}(x-b t-c).
    \label{eq:burgers_solution}
\end{equation}
Therefore, it is an ideal non-periodic benchmark to test the performance of the BSPF method. 

To solve Eq. \ref{eq:1d_burgers} numerically,  we rearrange the nonlinear advection term to the RHS, so that the spatial derivation can be applied entirely on the RHS, turning it into an initial value problem with the exact solution at $t=0$ as the initial condition: 
\begin{equation}
    \frac{\partial u}{\partial t} =RHS(u) =  \nu \mathcal{D}_{xx}u - u \mathcal{D}_xu, u(x, 0) = u_{exact}(x,0).
    \label{eq:1d_burgers_init}
\end{equation}
We choose the explicit Runge–Kutta method of order 5 (4) (RK45) as the time integration method and an exponential filter \citep{hou2007computing} for dealiasing. The time integration is carried out from 0 to 2 s. The problem is solved on the interval $[0,2\pi]$ using  the BSPF method ($p = 8, n = 32$, $N = 800$, uniform knots) for spatial discretizations. The initial condition is obtained from Eq. \ref{eq:burgers_solution} with the parameters of $a = 0.4$, $b =0.6$, $c = \pi$, $\nu = 0.01$, and the boundary conditions are enforced at $x=0$ and $2\pi$ with the analytical solution of Eq. \ref{eq:burgers_solution} in every timestep. 

\begin{figure}[h!]
\centering
\includegraphics[width=1\linewidth]{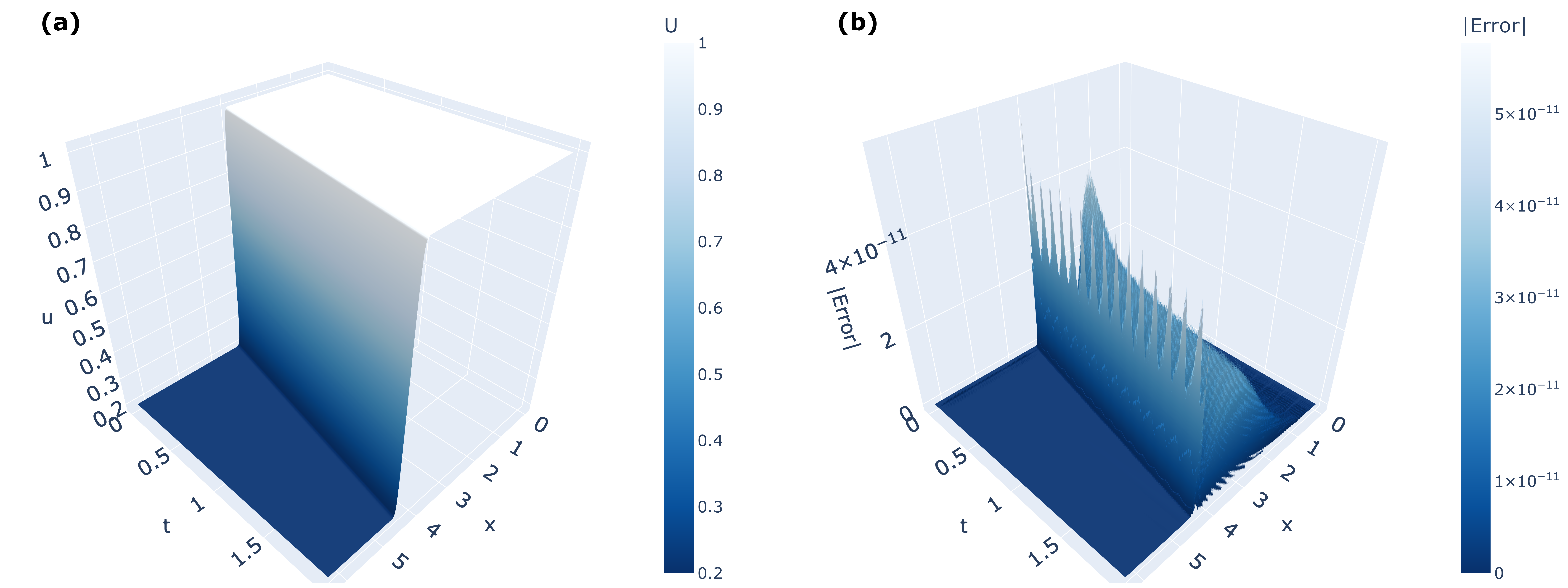}
\caption{\label{fig:burgers} Solving 1D Burgers' equation with the BSPF method: (a) the traveling wave solution; (b) distribution of the absolute error in the $x-t$ phase space. }
\end{figure}

Figure \ref{fig:burgers} shows the numerical solution of the BSPF method to Eq. \ref{eq:1d_burgers} in the $x-t$ plane (Figure \ref{fig:burgers}a) and the errors to the exact solutions (Figure \ref{fig:burgers}b). We can see that the solution of Eq. \ref{eq:burgers_solution} with the given initial condition Eq. \ref{eq:burgers_solution} is a traveling shock wave, and the wave front shape stays the same as it propagates. With $N = 800$, the largest absolute error to the exact solution is below $5\times10^{-11}$ and is observed at the location of the step due to the strong gradient.

\begin{figure}[h!]
\centering
\includegraphics[width=1\linewidth]{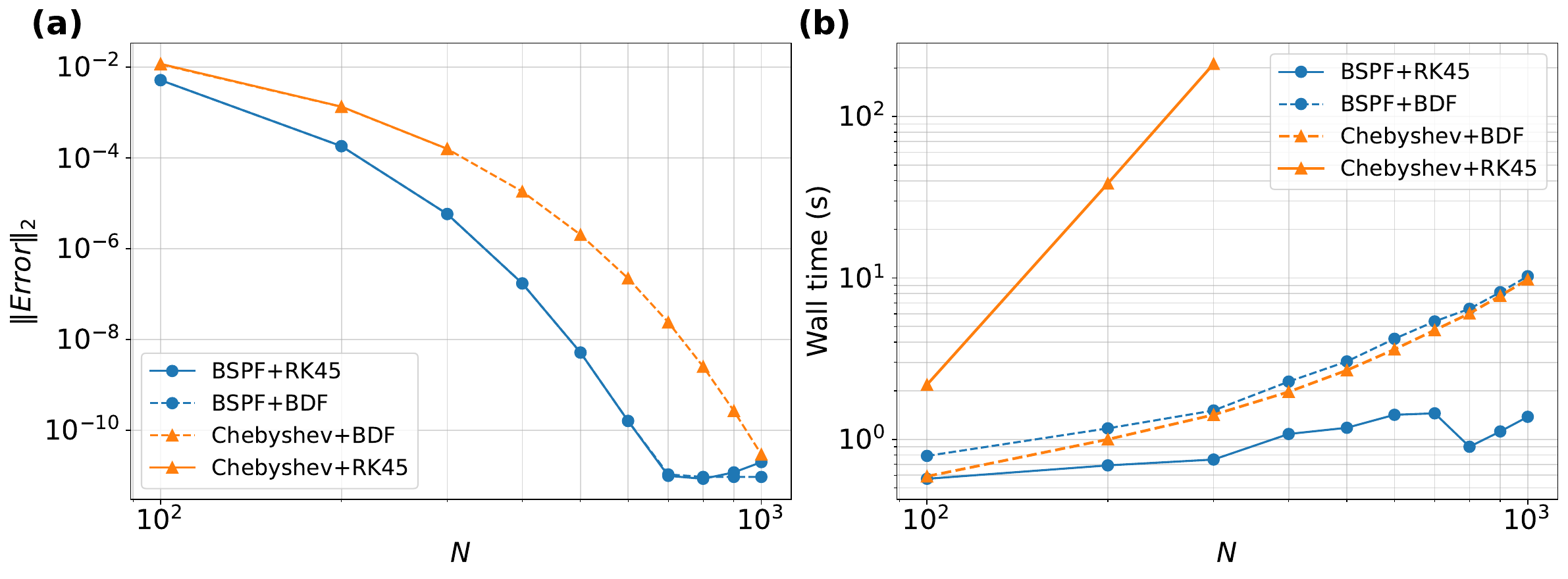}
\caption{\label{fig:burgers_convergence}Comparison of the Chebyshev method and the BSPF method in solving the 1D Burgers’ equation. (a) mesh convergence (b) time complexity.}
\end{figure}

Figure \ref{fig:burgers_convergence} further compares the performance of the BSPF and Chebyshev methods in solving the 1D Burgers' equation described by Eq. \ref{eq:1d_burgers} using the explicit RK45 and the implicit backward differentiation formula (BDF) method. From Figure  \ref{fig:burgers_convergence}a we can observe that the BSPF converges faster than the Chebyshev method. This is again due to the fact that the boundary-clustering Chebyshev node distribution is inefficient for capturing sharp jumps within the domain. Moreover, in Figure  \ref{fig:burgers_convergence}b we can see that the small grid size of the Chebyshev grid near the boundary impose a stiff stability constraint to the time integration scheme, requiring the usage of the implicit BDF scheme to solve it efficiently. The application of implicit schemes is more computationally expensive than the explicit RK45 schemes. As a result, despite the fact the efficiencies of differential operations for the BSPF and the Chebyshev methods  are similar, the Chebyshev method with the BDF scheme is significantly slower than the BSPF with the RK45 scheme, particularly for larger problem sizes.
\subsection{Two-Dimensional Shallow Water Equations}

In this section, we apply the BSPF method to solve the two-dimensional shallow water equations, a nonlinear PDE of great importance in applied mathematics, oceanography, and coastal engineering. This test case demonstrates the capability of the BSPF method to handle multidimensional problems. The 2D version of the BSPF method can be regarded as a “batched" application of the 1D BSPF operator in each spatial direction, and therefore inherits the same numerical behavior and accuracy as the 1D formulation.

The governing shallow water equations are formulated as \citep{imamura2006tsunami}:
\begin{equation}
\begin{aligned}
& \frac{\partial \eta}{\partial t}+\frac{\partial \mathcal{M}}{\partial x}+\frac{\partial \mathcal{N}}{\partial y}=0, \\
& \frac{\partial \mathcal{M}}{\partial t}+\frac{\partial}{\partial x}\left(\frac{\mathcal{M}^2}{\mathcal{H}}\right)+\frac{\partial}{\partial y}\left(\frac{\mathcal{M}\mathcal{N}}{\mathcal{H}}\right)+g \mathcal{H} \frac{\partial \eta}{\partial x}
  +\frac{g \alpha^2}{\mathcal{H}^{7/3}} \mathcal{M}\sqrt{\mathcal{M}^2+\mathcal{N}^2}=0, \\
& \frac{\partial \mathcal{N}}{\partial t}+\frac{\partial}{\partial x}\left(\frac{\mathcal{M}\mathcal{N}}{\mathcal{H}}\right)+\frac{\partial}{\partial y}\left(\frac{\mathcal{N}^2}{\mathcal{H}}\right)+g \mathcal{H} \frac{\partial \eta}{\partial y}
  +\frac{g \alpha^2}{\mathcal{H}^{7/3}} \mathcal{N}\sqrt{\mathcal{M}^2+\mathcal{N}^2}=0,
\end{aligned}
\end{equation}
where $\eta$ denotes the wave amplitude, $\mathcal{M}$ and $\mathcal{N}$ are the discharge fluxes in the $x$ and $y$ directions, respectively, and $\mathcal{H}$ is the total water column height defined by $\mathcal{H}=h+\eta$. where $h$ and $eta$ are given by Eqs. \ref{eq:h} and \ref{eq:eta}, respectively. The Manning roughness coefficient $\alpha$ is chosen as $\alpha=0.025$.

\begin{figure}[h!]
\centering
\includegraphics[width=1\linewidth]{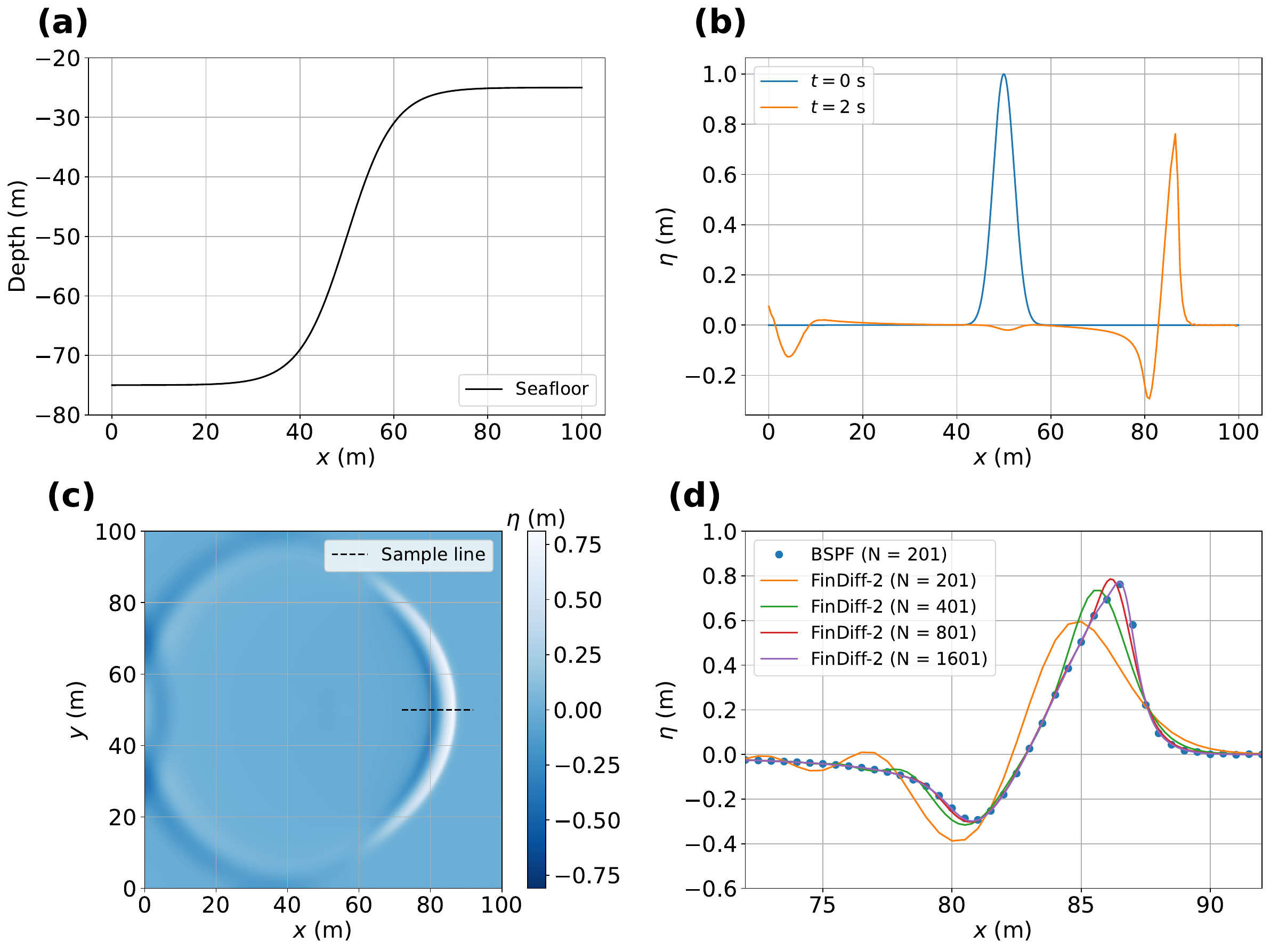}
\caption{\label{fig:tsunami} Simulation of the 2D shallow water equations using the BSPF method: (a) sea floor profile; (b) the surface elevation of y=50 m, the waveforms at t=0 (blue) and t=2 s (orange); (c) wave propagation at a later time (t = 2 s); (d) comparison of BSPF with finite difference (FD) results.}
\end{figure}

We consider a square computational domain of size $100~\text{m} \times 100~\text{m}$ with impenetrable wall boundary conditions. The BSPF discretization is configured with $\{p=8, n=32, N=201\times 201\}$ and uniform knot spacing. Time integration is performed using a fourth-order explicit Runge--Kutta scheme with a fixed step size of $0.001~\text{s}$ (Figure~\ref{fig:tsunami}).

The sea floor is defined by a hyperbolic tangent profile to mimic a continental shelf  (Figure~\ref{fig:tsunami}a):
\begin{equation}
h(x) = 50 - 25 \tanh\!\left(\frac{x-50}{10}\right).
\label{eq:h}
\end{equation}
The initial free-surface elevation is prescribed as a Gaussian pulse located at the domain center:
\begin{equation}
\eta(x,y,0) = \exp\!\left[-\frac{(x-50)^2}{10} - \frac{(y-50)^2}{10}\right].
\label{eq:eta}
\end{equation}

At $t=0$ s, the initial Gaussian water column collapses and excites the propagation of shallow water waves. Due to the decreasing water depth in the positive $x$-direction, the nonlinear dynamics of the shallow water system generate a sharp wave front in this direction, which presents a significant challenge for numerical schemes (Figure~\ref{fig:tsunami}b \& Figure~\ref{fig:tsunami}c).

To validate the BSPF results, we compare the wave profile along $y=0$ with the solutions obtained using a classical second-order finite difference (FD) solver \cite{koehn2025shallowwater}, with grid resolutions ranging from $N=201$ to $N=1601$. As illustrated in Figure~\ref{fig:tsunami}d, the BSPF solution at $N=201$ closely matches the FD solution at $N=1601$. In contrast, FD solutions with lower grid resolution fail to resolve the sharp wave front, thereby demonstrating that the BSPF method achieves high-order numerical accuracy while requiring significantly fewer grid points.

\section{Discussion and Conclusions}

In this work, we have proposed and demonstrated the B-spline-periodized Fourier (BSPF) pseudo-spectral method as a general framework for solving non-periodic problems with spectral-like accuracy. By combining a boundary-constrained B-spline approximation with Fourier-based residual correction, the BSPF method enables smooth and high-order periodization of non-periodic functions without the need to enlarge the computational domain. This approach allows us to retain the efficiency of FFT-based spectral algorithms while mitigating the Gibbs phenomenon at the boundaries.

Our theoretical analysis has shown that the accuracy of the BSPF method is determined by two main factors: the order of the B-spline approximation and the accuracy of the boundary derivative stencils. The method achieves near-spectral convergence when the errors is dominant by the interior part, while the finite-difference approximation of the boundary derivatives leads to algebraic convergence consistent with the stencil order. A series of numerical tests on differentiation tasks confirmed the theoretical error bounds and demonstrated the robustness of the BSPF method under challenging conditions such as noisy and oscillatory functions. Comparisons with Chebyshev pseudo-spectral and high-order finite-difference methods reveal that BSPF offers superior accuracy on interior-oscillatory data, while retaining comparable computational complexity $\mathcal{O}(N \log N)$ due to FFT dominance. Moreover, the ability to analytically map BSPF to non-uniform grids further enhances its adaptability to problems where local refinement is required.

We have also applied the BSPF method to nonlinear PDEs, including the one-dimensional Burgers' equation and the two-dimensional shallow water equations. In both cases, the BSPF method successfully captured sharp gradients and nonlinear wave dynamics with high accuracy. The shallow water test, in particular, highlighted the advantage of BSPF over lower-order methods, achieving results comparable to high-resolution simulations at significantly lower grid numbers.

Overall, the BSPF method provides the following advantages:
\begin{enumerate}[leftmargin=*]
    \item High-order accuracy for non-periodic problems while maintaining FFT-based efficiency.
    \item Straightforward enforcement of boundary conditions via boundary-constrained B-spline fitting.
    \item Flexibility to operate on both uniform and analytically mapped non-uniform grids.
    \item Robustness against noisy data and superior performance compared to existing spectral and finite-difference approaches in capturing interior-oscillatory phenomena.
\end{enumerate}

Future research directions include the exploration of optimal knot distributions for improved stability and accuracy, the extension of BSPF to more complex geometries and boundary conditions, and its application to large-scale multidimensional problems in fluid dynamics, oceanography, and plasma physics. In addition, the integration of adaptive refinement strategies and parallel implementations will further enhance the efficiency and applicability of the BSPF method in high-performance computing environments.

In conclusion, the BSPF pseudo-spectral method provides a powerful and flexible tool for solving non-periodic PDE problems with spectral-like accuracy and efficiency, bridging the gap between Fourier and polynomial-based spectral methods.

\section*{CRediT authorship contribution statement}
\textbf{Dongan Li}: Conceptualization, Methodology, Software, Validation, Formal analysis, Investigation, Writing - Original Draft, Visualization. \textbf{Mou Lin}: Methodology, Software, Validation, Formal analysis, Investigation, Data Curation, Writing - Review \& Editing , Visualization. \textbf{Shunxiang Cao}: Writing - Review \& Editing. \textbf{Shenli Chen}: Investigation,  Writing - Review \& Editing, Resources, Supervision, Project administration, Funding acquisition.

\section*{Data availability}
The data sets and programs used in this study are publicly accessible at \\
https://github.com/moulin1024/bspf. 

\section*{Declaration of competing interest}
The authors declare that there is no conflict of interest.

\section*{Acknowledgment}

This study is supported by funds from the National Key R\&D Program of China (2024YFC3013201), Shenzhen Science and Technology Innovation Committee (KCXFZ20240903093900002), and Guangdong Basic and Applied Basic Research Foundation (2022B1515130006).

 \bibliographystyle{elsarticle-num} 
 \bibliography{ref}






\end{document}